\magnification=\magstep1    
\def\m@th{\mathsurround=0pt}
\def\recase#1{\phantom\{ \vcenter{\normalbaselines\m@th
     \ialign{##\hfil & \quad ## \crcr#1\crcr}}\Biggr\} }
\def\recasebig#1{\left. \vcenter{\normalbaselines\m@th
     \ialign{##\hfil & \quad ## \crcr#1\crcr}}\right\} }

\def \R{{\bf R}}
\def \N{{\bf N}}
\def \Z{{\bf Z}}

\def \cD{{\cal D}}

\def \cH{{\cal H}}

\def \cL{{\cal L}}

\def \cF{{\cal F}}

\def \unp3{{\rm unps3}}

\def \card{{\rm card}}

\hyphenation{argu-ment}
\def\ref{\bigskip\hangindent=25pt\hangafter=1\noindent}
\def\refs{\medskip\hangindent=25pt\hangafter=1\noindent}
\overfullrule= 0 pt 
\parskip =0pt

\vglue.5in
 \centerline{\bf A STRICTLY STATIONARY, $N$-TUPLEWISE
INDEPENDENT} \centerline{\bf COUNTEREXAMPLE TO THE CENTRAL LIMIT
THEOREM}
\bigskip
\bigskip

$$\vbox{\halign{& # \hfil \quad & #\hfil \cr & Richard C.\ Bradley &
Alexander R.\ Pruss \cr &Department of Mathematics
$\qquad\qquad\quad$& Department of Philosophy \cr &Indiana
University & Baylor University \cr &Rawles Hall &One Bear Place
\#97273 \cr &Bloomington, Indiana 47405 & Waco, Texas 76798-7273 \cr
&USA & USA\cr} }$$

\bigskip
\vskip .5in

\par{\bf Abstract.} For an arbitrary integer $N\ge 2$, this paper
gives a construction of a strictly stationary, $N$-tuplewise
independent sequence of (nondegenerate) bounded random variables
such that the Central Limit Theorem fails to hold.  The sequence is
in part an adaptation of a nonstationary example with similar
properties constructed by one of the authors (ARP) in a paper
published in 1998.
 \medskip
\bigskip
\vskip 2.5in

\par {\it AMS} 2000 {\it Mathematics Subject Classification.} 60G10, 60F05.

 \vskip .2truein \vfill\eject


\noindent {\bf 1. Introduction and main result} \hfil\break

Suppose $X := (X_k, k \in \Z)$ is a sequence of random variables on
a probability space $(\Omega, {\cal F}, P)$. This sequence $X$ is
said to be ``strictly stationary'' if for all choices of integers
$j$ and $\ell$ and nonnegative integer $m$, the random vectors
$(X_j, X_{j+1}, \dots, X_{j+m})$ and $(X_\ell, X_{\ell+1}, \dots,
X_{\ell+m})$ have the same distribution.   For a given integer $N
\geq 2$, the sequence $X$ (stationary or not) is said to satisfy
``$N$-tuplewise independence'' if for every choice of $N$ distinct
integers $k(1),\ k(2),\ \dots,\ k(N)$, the random variables
$X_{k(1)},\ X_{k(2)},\ \dots,\ X_{k(N)}$ are independent. For $N=2$
(resp.\ $N=3$), the word ``$N$-tuplewise'' is also
expressed as ``pairwise'' (resp.\ ``triplewise''). 

  Etemadi [5] proved a strong law of large numbers
for sequences of pairwise independent, identically distributed
random variables with finite absolute first moment. Janson [7]
showed with several classes of counterexamples that for strictly
stationary sequences of pairwise independent, nondegenerate,
square-integrable random variables, the Central Limit Theorem
(henceforth abbreviated CLT) need not hold. Subsequently, Bradley
[2, Theorem 1] constructed another such counterexample, a 3-state
one that has the additional property of satisfying the absolute
regularity (weak Bernoulli) condition.  For an arbitrary fixed
integer $N \geq 3$, Pruss [9] constructed a (not strictly
stationary) sequence of bounded, nondegenerate, $N$-tuplewise
independent, identically distributed random variables for which the
CLT fails to hold. In that paper, Pruss left open the question
whether, for any integer $N \geq 3$, a strictly stationary
counterexample exists.     For $N = 3$, Bradley [3, Theorem 1]
answered that question by showing that the counterexample in [2,
Theorem 1] alluded to above is in fact triplewise independent.

  In a similar spirit, for an arbitrary integer
$N \geq 2$, Flaminio [6] constructed a strictly stationary,
finite-state, $N$-tuplewise independent random sequence $X := (X_k,
k \in \Z)$ which also has zero entropy and is mixing (in the
ergodic-theoretic sense). That paper explicitly left open the
question of whether those examples satisfy the CLT. 

  In this paper, we shall answer affirmatively the question in
[9], by constructing for an arbitrary fixed integer $N \geq 2$ a
strictly stationary, $N$-tuplewise independent sequence of bounded,
nondegenerate random variables such that the CLT fails to hold.  The
construction will be in part an adaptation of the (not strictly
stationary) counterexample in [9].

  Before the main result is stated, a few notations are needed:

  The Borel $\sigma$-field on the real number line ${\R}$
will be denoted ${\cal R}$.

  Convergence in distribution will be denoted by $\Rightarrow$.

  The set of positive integers will be denoted by ${\N}$.   For a given sequence $X :=
(X_k, k \in {\Z})$ of random variables, the partial sums will be
denoted for $n \in {\N}$ by
$$ S_n := S(X,n) := X_1 + X_2 + \dots + X_n.  \eqno (1.1)$$

  Here is our main result: \hfil\break

\noindent {\bf Theorem 1.1.}\ \ {\sl Suppose $N$ is an integer such
that $N \geq 2$.  Then there exists a strictly stationary sequence
$X := (X_k, k \in {\Z})$ of random variables such that the following
statements hold:

  (A) The random variable $X_0$ is uniformly distributed
on the interval $[-3^{1/2}, 3^{1/2}]$ (and hence $EX_0 = 0$ and
$EX_0^2 = 1$).

  (B) For every choice of $N$ distinct integers
$k(1),\ k(2),\ \dots,\ k(N)$, the random variables $X_{k(1)},\
X_{k(2)},\ \dots,\ X_{k(N)}$ are independent.

  (C) The random variables $|X_k|,\ k \in {\Z}$ are
independent (and identically distributed).

  (D) For every infinite set $Q \subset {\N}$,
there exist an infinite set $T \subset Q$ and a nondegenerate,
non-normal probability measure $\mu$ on $({\R}, {\cal R})$ such that
$ S_n/n^{1/2} \Rightarrow \mu\ \ {\rm as}\ \ n \to \infty,\ n \in
T$.}
\bigskip

   Here are some comments on Theorem 1.1:

  By property (A), the ``natural normalization'' to consider
for the central limit question for the partial sums of this sequence
is $S_n/n^{1/2}$.

  Property (B) is of course $N$-tuplewise independence.

  By property (D) and the Theorem of Types (see e.g.\
[1, p.\ 193, Theorem 14.2]), there do not exist constants $a_n > 0$
and $b_n \in {\R}$ for $n \in {\N}$ such that $a_n S_n + b_n
\Rightarrow N(0,1)$, even along a subsequence of the positive
integers.

  Also by property (D) and an elementary argument, there do not
exist constants $b_n \in {\R}$ for $n \in {\N}$ such that
$n^{-1/2}S_n + b_n \to 0$ in probability, even along a subsequence
of the positive integers.

  In property (D), the probability measure $\mu$ may depend on
the set $Q$. \hfil\break

\noindent {\bf Remark 1.2.}\ \ With essentially the same
construction, one obtains an analog of Theorem 1.1 with property (A)
replaced by the following one:  (A$'$) The random variable $X_0$ has
the $N(0,1)$ distribution. \hfil\break

  We did not investigate the ergodic-theoretic properties of
the sequence $X$ in Theorem 1.1, nor did we try to ascertain the
particular class of probability measures $\mu$ that can arise in
statement (D) there. Bradley [4] gives a construction of a
(nondegenerate, two-state) strictly stationary,  5-tuplewise
independent random sequence which fails to satisfy the CLT (instead,
it satisfies property (D) in Theorem 1.1), and also has the extra
properties of being ``causal'' (for an appropriate use of that term)
and therefore ``Bernoulli'' (i.e.\ isomorphic to a Bernoulli shift)
and also having a trivial double tail $\sigma$-field.  There does
not seem to be an obvious way to build those extra properties into
the construction here for Theorem 1.1. \hfil\break

  The proof of Theorem 1.1 will be given in sections 2 and 3.
It will be rather intricate and will require a fair amount of
notation.  Refer to properties (A) and (C) in Theorem 1.1.  The
construction will basically involve taking a sequence of
independent, identically distributed random variables uniformly
distributed on the interval $[-3^{1/2}, 3^{1/2}]$, and changing the
signs of the variables in such a way as to introduce a dependence
that preserves $N$-tuplewise independence but is incompatible with
the CLT. \hfil\break

In essence it involves the conversion of the (not strictly
stationary) counterexample in [9] to one that is strictly
stationary.  The main ideas for that conversion were outlined in an
e-mail message by one of the authors (Pruss [10]) to the other
author (RCB), and are developed in section 3 here. Section 2 gives
some vital ``preliminary'' information; much of it is taken or
adapted from [9], but will be given again here in detail because of
considerable extra complications in our context. The use (in both
sections 2 and 3) of higher order moments to establish property (D)
in Theorem 1.1, is adapted from the analogous use of 6th moments for
the same purpose in (an earlier, 2006 version of) the preprint [4].
Underlying all this is the repeated creation of ``big'' collections
of $N$-tuplewise independent random variables from ``smaller'' ones;
such procedures are well known in the theory of error-correcting
codes (see e.g.\ [8]).
\bigskip

\noindent{\bf 2. Part 1 of proof of Theorem 1.1}
\medskip

Sections 2 and 3 together will give the proof of Theorem 1.1. Both
sections will be divided into several ``steps,'' including a
``definition,'' some ``lemmas,'' etc. Throughout this proof, the
setting is a probability space $(\Omega,\cF,P)$, ``enlarged'' as
necessary to accommodate all random variables defined in this proof.
\medskip

\noindent{\bf Step 2.1.} Let $L$ be an arbitrary fixed integer such
that
$$\hbox{$L$ is even and $L\ge 6$.} \eqno{(2.1)}$$
To prove Theorem 1.1, it suffices to construct a strictly stationary,
$(L-1)$--tuplewise independent random sequence $X:= (X_k$, $k\in\Z)$
that also satisfies properties  (A), (C), and (D) in Theorem 1.1.
That will be the goal of sections 2 and 3.
\medskip

\noindent {\bf Step 2.2.} The following notations and conventions
will be used:

(a) Refer to (2.1). For $n\in \{0,1,2,\dots\}$, when the term $L^n$
appears in a subscript or exponent, it will be written as $L\uparrow
n$ for typographical convenience.

(b) Suppose $n\in\N$. A vector $x\in\R^n$ will often be represented
as $x :=(x_0,x_1,\dots,x_{n-1})$ (instead of $(x_1,x_2,\dots,x_n))$,
for ``bookkeeping'' convenience. For a given $x:= (x_0,x_1,\dots,$
$x_{n-1})\in\R^n$, define the two real numbers
$$\hbox{sum}\,\,x := \sum^{n-1}_{i=0}x_i\quad \hbox{and}\quad
\hbox{prod}\,\,x :=\prod^{n-1}_{i=0}x_i. \eqno{(2.2)}$$ Notations of
the form $f((x_0,x_1,\dots,x_{n-1}))$ will be written simply as
$f(x_0$, $x_1,$ $\dots,x_{n-1})$.

(c) Sometimes the coordinates of a vector will be permuted. If
$n\in\N$, $x:= (x_0,x_1,\dots,$ $x_{n-1})\in\R^n$, and $\sigma$ is a
permutation of the set $\{0,1,\dots,n-1\}$, then define the vector $
x_\sigma \in\R^n$ by $x_\sigma :=
((x_\sigma)_0,(X_\sigma)_1,\dots,(X_\sigma)_{n-1}) :=
(x_{\sigma(0)},x_{\sigma(1)},\dots,x_{\sigma(n-1)})$.

(d) If $(x_k$, $k\in\Z)\in\R^\Z$, and $a\le b$ are integers, then
the vector $(x_a,x_{a+1},\dots,x_b)$ will also be denoted $(x_k:
a\le k\le b)$.

(e) A ``measure'' on the space $\R$ or $\R^n$ $(n\in\N)$ or $\R^\Z$
will always mean a measure on the Borel $\sigma$-field on that
space.

(f) If $\eta$ is a random variable/vector/sequence, then the
distribution of $\eta$ on the appropriate space (such as $\R$,
$\R^n$, or $\R^\Z$) will be denoted $\cL(\eta)$. If also $F$ is an
event such that $P(F)>0$, then $\cL(\eta\mid F)$ will denote the
conditional distribution of $\eta$, given~$F$.

(g) For a given $n\in\N$, an ``$\R^n$--valued random vector'' is
simply a random vector with $n$ coordinates. Of course, if $n,m\in
\N$, $\eta_1$ and $\eta_2$ are $\R^n$--valued random vectors,
$\cL(\eta_1) = \cL(\eta_2)$ (on $\R^n$), and $f:\R^n\to\R^m$ is a
Borel function, then $\cL(f(\eta_1)) = \cL(f(\eta_2))$ (on $\R^m$).

(h) Several notations will be defined here. Suppose $Y:= (Y_i$,
$i\in I)$ is a family of random variables, where $I$ is a nonempty
(possibly infinite) index set. The $\sigma$--field of events
generated by this family will be denoted $\sigma(Y)$ or
$\sigma(Y_i$, $i\in I)$. To avoid any confusion, this family $Y$ is
said to satisfy ``$(L-1)$--tuplewise independence'' if either
(i)~card $I=1$, or (ii)~card $I\ge 2$ and for every set $S\subset I$
such that $2\le \card$ $S\le L-1$, the random variables $Y_i$, $i\in
S$ are independent. (The point there is to formally include the case
card $I=1$ in that terminology.) Here and below, ``card'' means
cardinality. In the case of a random vector $Y:= (Y_0,Y_1,\dots,
Y_{n-1})$, where $n\in\N$, the phrase ``$Y$ satisfies
$(L-1)$--tuplewise independence'' simply means that the family of it
coordinates $(Y_i$, $i\in \{0,1,\dots,n-1\})$ satisfies
$(L-1)$--tuplewise independence.
\medskip

\noindent{\bf Step 2.3.} Here a few useful functions on $\R^n$ will
be defined.
\medskip

\noindent{\bf Step 2.3(A).} For each $n\in\N$, define the function
$\varphi_n :\R^n\to \R^n$ as follows: For $x\in\R^n$ (see (2.2)),
$$\varphi_n(x) := \cases{
x & if sum $x>0$ \cr 0_n & if sum $x=0$ \cr -x & if sum $x<0$ \cr}
\eqno{(2.3)}$$ where $0_n$ denotes the origin in $\R^n$.

Thus for example $\varphi_3(5,-7,4) = \varphi_3(-5,7,-4) = (5,-7,4)$.
\medskip

{\it Remark 1.} For each $n\in\N$ and each $x\in\R^n$, sum
$\varphi_n(x) = |\hbox{sum $x$}|$.
\medskip

{\it Remark 2.} Suppose $n\in\N$, $x:=
(x_0,x_1,\dots,x_{n-1})\in\R^n$, $y := \varphi_n(x) :=
(y_0,y_1,\dots,$ $y_{n-1})$, and $\sigma$ is a permutation of the
set $\{0,1,\dots,n-1\}$. Recall from step 2.2(c) the notations
$x_\sigma := (x_{\sigma(0)},x_{\sigma(1)},\dots,x_{\sigma(n-1)})$
and $y_\sigma :=
(y_{\sigma(0)},y_{\sigma(1)},\dots,y_{\sigma(n-1)})$. Then
 $y_\sigma= \varphi_n(x_\sigma)$. This holds by a careful trivial
 argument, using the fact that sum $x_\sigma =$ sum~$x$. Thus for
 $0\le i\le n-1$ (see step 2.2(c) again), $(\varphi_n(x_\sigma))_i =
 (y_\sigma)_i = y_{\sigma(i)} = (\varphi_n(x))_{\sigma(i)}$.
\medskip

{\it Remark 3.} If $n,x$, and $y$ are as in Remark~2 above (with $y=
\varphi_n(x))$, and also sum~$x\not= 0$, then $|x_i| = |y_i|$ for
each $i\in \{0,1,\dots,n-1\}$.
\medskip

{\it Remark 4.} Suppose $n\in\N$, $Y$ is an $\R^n$--valued random
vector (recall step 2.2(g)) such that $\cL(-Y) = \cL(Y)$ and
$\cL(Y)$ is absolutely continuous with respect to Lebesgue measure
on $\R^n$, and $V$ is a random variable independent of $Y$ such that
$P(V=-1) = P(V=1) =1/2$. Then $P(\hbox{sum}\,\,Y=0) =0$; and for any
Borel set $B\subset \{x\in\R^n: \hbox{sum}\,\,x>0\}$,
$$P(Y\in B) = P(-Y\in B) = (1/2) \cdot P(\varphi_n(Y)\in B);$$
and hence by a trivial argument, $\cL(V\varphi_n(Y)) = \cL(Y)$.
\medskip

\noindent{\bf Step 2.3(B).} Refer to (2.1). Occasionally, ``big''
vectors will be created by the splicing together of $L$ ``smaller''
ones. Suppose $n$ is a positive integer, and for each
$\ell\in\{0,1,\dots,L-1\}$, $x^{(\ell)}:=
(x^{(\ell)}_0,x^{(\ell)}_1,\dots,x^{(\ell)}_{n-1})\in \R^n$. Then the
notation
$$\left\langle x^{(0)} \mid x^{(1)} \mid x^{(2)} \mid \cdots \mid
x^{(L-1)}\right\rangle \eqno{(2.4)}$$
means the vector $y := (y_0,y_1,y_2,\dots,y_{Ln-1}) \in \R^{Ln}$ such
that for each $\ell\in \{0,1,\dots,L-1\}$, $(y_{\ell n}, y_{\ell
n+1},\dots,y_{\ell n+n-1}) = x^{(\ell)}$.
\medskip

\noindent{\bf Step 2.3(C).} For each $n\in\N$ and each $j\in
\{0,1,\dots,L-1\}$, define the function $\psi_{n,j}: \R^{Ln} \to
\R^{Ln}$ as follows: Suppose $y\in \R^{Ln}$. Represent $y := \langle
x^{(0)} \mid x^{(1)} \mid \cdots \mid x^{(L-1)}\rangle$ as in (2.4),
where for each $\ell \in \{0,1,\dots,L-1\}$, $x^{(\ell)}\in\R^n$.
Then (see (2.2)) define
$$\psi_{n,j}(y) := \cases{ y &if $\prod^{L-1}_{\ell=0}(\hbox{sum}\,
x^{(\ell)}) \le 0$ \cr \langle x^{(0)}|\cdots
|x^{(j-1)}|-x^{(j)}|x^{(j+1)}|\cdots|x^{(L-1)}\rangle &if
$\prod^{L-1}_{\ell=0}(\hbox{sum}\, x^{(\ell)})>0$. \cr}
\eqno{(2.5)}$$ That is, if $\prod^{L-1}_{\ell=0} ($sum
$x^{(\ell)})>0$, then to define $\psi_{n,j}(y)$, one just multiplies
by $-1$ the coordinates of $x^{(j)}$ (the $(j+1)^{st}$
``$n$--coordinate piece'' of $y$); all other coordinates of $y$ are
left alone.
\medskip

{\it Preview.}  We shall show later on that if for each $\ell\in
\{0,1,\dots,L-1\}$, $Y^{(\ell)}$ is an $\R^n$--valued random vector
satisfying $(L-1)$--tuplewise independence (see (2.1)) and certain
other conditions, and these random vectors $Y^{(\ell)}$, $\ell\in
\{0,1,\dots,L-1\}$ are independent of each other, and one defines
the $\R^{Ln}$--valued random vector $Y:= \langle Y^{(0)}\mid
Y^{(1)}\mid \cdots \mid Y^{(L-1)}\rangle$, then for any $j\in
\{0,1,\dots,L-1\}$, the $\R^{Ln}$--valued random vector
$\psi_{n,j}(Y)$ satisfies $(L-1)$--tuplewise independence. Iterative
use of this construction will generate the final sequence $X:=
(X_k$, $k\in\Z)$ for Theorem~1.1.
\medskip

\noindent{\bf Step 2.4.} Here several special probability measures on
$\R$ or $\R^n$ or $\R^\Z$ will be defined.
\medskip

\noindent{\bf Step 2.4(A).} The notation $\lambda_{\unp3}$ will
refer to the uniform distribution on the interval
$[-3^{1/2},3^{1/2}]$, regarded as a probability measure on $\R$.
(The subscript ``unps3'' stands for {\bf u}niform on the interval
from {\bf n}egative to {\bf p}ositive {\bf s}quare root of {\bf 3}.)
For any positive integer $m$, the $m$--fold product measure
$\lambda_{\unp3} \times \lambda_{\unp3} \times \dots \times
\lambda_{\unp3}$ on $\R^m$ will be denoted $\lambda^{[m]}_{\unp3}$.
\medskip

{\it Remark.} If $U$ is a random variable such that $\cL(U) =
\lambda_{\unp3}$ (see step 2.2(f)), then of course $EU^n=0$ for odd
$n\in \N$ and $EU^n = (n+1)^{-1} \cdot 3^{n/2}$ for even $n\in \N$.
Hence, if also $Z$ is a $N(0,1)$ random variable, then (see e.g.\
[1, p. 275, eq.\ (21.7)]) for every $n\in\N$, $EZ^n \ge EU^n \ge 0$
(with equality for odd $n$ and with $EZ^2 = EU^2 =1$).
\medskip

\noindent{\bf Step 2.4(B).} Suppose $n\in\N$, and $\mu$ is a
probability measure on $\R^n$. Then the notation $\cD(n,\mu)$ means
the distribution on $\R^\Z$ of a random sequence $Y:= (Y_k$,
$k\in\Z)$ such that (i)~for each $u\in\Z$, the random vector
$\zeta^{(u)}:= (Y_{nu}, Y_{nu+1},Y_{nu+2},\dots,Y_{nu+n-1})$
satisfies $\cL(\zeta^{(u)}) = \mu$, and (ii)~these random vectors
$\zeta^{(u)}$, $u\in\Z$ are independent of each other.
\medskip

\noindent{\bf Step 2.4(C).} Define (see (2.1) and (2.2)) the set
$\Upsilon := \{x\in \{-1,1\}^L: \hbox{prod}\thinspace x = -1\}$.
Then card~$\Upsilon = 2^{L-1}$. Let $\nu$ denote the uniform
distribution on $\Upsilon$
--- that is, the probability measure $\nu$ on $\R^L$ such that $\nu(\{x\})
= 1/2^{L-1}$ for each $x \in \Upsilon$.
\medskip

\noindent{\bf Step 2.4(D).} Suppose $n\in\N$, and $\mu$ is a
probability measure on $\R^n$. Then define the probability measure
$\theta(\mu)$ on $\R^{Ln}$ by $\theta(\mu) := \cL(Y)$ where the
random vector $Y:= (Y_0,Y_1,\dots, Y_{Ln-1})$ is as defined below:

Let $W^{(0)},W^{(1)},\dots,W^{(L-1)}$ each be an $\R^n$--valued
random vector, such that $\cL(W^{(\ell)}) =\mu$ for each $\ell \in
\{0,1,\dots,L-1\}$, and such that these random vectors $W^{(\ell)}$,
$\ell\in \{0,1,\dots,L-1\}$ are independent of each other. Let $V:=
(V_0,V_1,\dots, V_{L-1})$ be a $\Upsilon$--valued random vector such
that $\cL(V) = \nu$ (see step 2.4(C) above), with $V$ being
independent of the family $(W^{(\ell)}$,
$\ell\in\{0,1,\dots,L-1\})$. Referring to step 2.3(A)(B), define the
$\R^{Ln}$--valued random vector $Y$ by
$$Y := \left\langle V_0\varphi_n(W^{(0)})\mid V_1\varphi_n(W^{(1)})
\mid V_2 \varphi_n(W^{(2)})\mid \cdots \mid V_{L-1}\varphi_n
(W^{(L-1)})\right\rangle. \eqno{(2.6)}$$

\noindent{\bf Remark 2.5.} Suppose $V:= (V_0,V_1,\dots,V_{L-1})$ is
a random vector such that $\cL(V) = \nu$ (see (2.1) and step
2.4(C)). Then by elementary arguments, the following statements
hold:

(i) For each $k\in \{0,1,\dots,L-1\}$, $P(V_k =-1) = P(V_k=1) = 1/2$.

(ii) $V$ satisfies $(L-1)$--tuplewise independence.

(iii) $\cL(-V) = \nu$.

(iv) For every permutation $\sigma$ on the set $\{0,1,\dots,L-1\}$,
the random vector $V_\sigma :=
(V_{\sigma(0)},V_{\sigma(1)},\dots,V_{\sigma(L-1)})$ (see step
2.2(c)) satisfies $\cL(V_\sigma) = \nu$.

(v) prod $V=-1$ a.s.\ (see (2.2)). Hence for any $j\in
\{0,1,\dots,L-1\}$, \break \noindent$V_j =-
\prod_{\ell\in\{0,1,\dots,L-1\}-\{j\}}V_\ell$ a.s.
\medskip

\noindent{\bf Definition 2.6.} Suppose $n\in\N$, and $\mu$ is a
probability measure on $\R^n$. Then $\mu$ is said to satisfy
``Condition $\cH(n)$'' if the following seven statements (a)--(g)
hold, where $W:= (W_0,W_1,\dots,W_{n-1})$ is a random vector such
that $\cL(W) = \mu$:

(a) $\mu$ is absolutely continuous with respect to Lebesgue measure
on~$\R^n$.

(b) For each $k\in \{0,1,\dots,n-1\}$, $\cL(W_k) = \lambda_{\unp3}$
(see step 2.4(A)).

(c) $W$ satisfies $(L-1)$--tuplewise independence.

(d) The random variables $|W_0|,|W_1|,\dots,|W_{n-1}|$ are
independent (if $n\ge 2$).

(e) $\cL(-W) =\mu$.

(f) For any $j\in \{0,1,\dots,n-1\}$, there exists a permutation
$\sigma$ on the set $\{0,1,\dots,n-1\}$ such that (i)~$\sigma(0)=j$
and (ii)~$\cL(W_\sigma)=\mu$ (where $W_\sigma:= (W_{\sigma(0)},
W_{\sigma(1)},\dots,W_{\sigma(n-1)})$ --- see step 2.2(c)).

(g) Either (i)~$n<L$, or (ii)~$n\ge L$, and for every set $S\subset
\{0,1,\dots,n-1\}$ such that card $S=L$, one has that $E(\prod_{i\in
S} W_i)\le 0$.
\medskip

\noindent{\bf Lemma 2.7.} {\sl Suppose $n\in\N$, and $W:=
(W_0,W_1,\dots,W_{n-1})$ is a random vector whose distribution
$\cL(W)$ on $\R^n$ satisfies Condition $\cH(n)$. Referring to (2.3),
define the random vector $Y:= \varphi_n(W):=
(Y_0,Y_1,\dots,Y_{n-1})$. Then the following three statements hold
(see (2.2)):

(i) $E(\hbox{\rm sum}\, Y) = E|\hbox{\rm sum}\, W|\ge (1/2)n^{1/2}$.

(ii) $\cL(Y_0) = \cL(Y_1) = \cdots =\cL(Y_{n-1})$.

(iii) For any $k\in \{0,1,\dots,n-1\}$, $EY_k = (1/n) \cdot
E|\hbox{\rm sum}\, W|$.}
\medskip

\noindent{\bf Proof.} Let us first prove (i).  Trivially sum $Y
=|$sum $W|$ by (2.3). We just need to prove the (``latter'')
inequality in (i).

Refer to (2.1) and Definition 2.6(b)(c). By simple calculations,
including the argument in [1, p.\ 85, proof of Theorem 6.1] (which
uses only 4--tuplewise independence), one has that for each $k\in
\{0,1,\dots,n-1\}$, $\cL(W_k) = \lambda_{\unp3}$, $EW_k =0$,
$EW^2_k=1$, and $EW^4_k =9/5$, and hence $E($sum $W)^2 =n$ and
$$
E(\hbox{sum }\, W)^4 = n\cdot EW^4_0 +3n(n-1) \cdot (EW^2_0)^2 <3n^2.
$$
Hence by H\"older's inequality,
$$\eqalign{
n= E(\hbox{sum }\, W)^2 &\le \left\| |\hbox{sum }\,
W|^{2/3}\right\|_{3/2} \cdot \left\| |\hbox{sum }\,
W|^{4/3}\right\|_3 \cr &= [E|\hbox{sum }\, W|\,]^{2/3} \cdot
\left[E(\hbox{sum }\, W)^4 \right]^{1/3} \cr &\le [E|\hbox{sum }\,
W|\,]^{2/3} \cdot \left(3n^2\right)^{1/3}. \cr}$$ Hence $n^{3/2} \le
E|\hbox{sum }\, W|\cdot (3n^2)^{1/2}$. Hence the (``latter'')
inequality in (i) holds.

Now let us prove (ii). Suppose $j\in \{0,1,\dots,n-1\}$. Referring to
Definition 2.6(f), let $\sigma$ be a permutation of the set
$\{0,1,\dots,n-1\}$ such that $\sigma(0) =j$ and $\cL(W_\sigma) =
\cL(W)$. Then $\cL(\varphi_n(W_\sigma)) = \cL(\varphi_n(W)) =
\cL(Y)$. Also, $Y_\sigma = \varphi_n (W_\sigma)$ by Remark~2 in step
2.3(A), and hence $\cL(Y_\sigma) = \cL(Y)$. Hence $\cL(Y_j) =
\cL(Y_{\sigma(0)}) = \cL(Y_0)$. Since $j\in \{0,1,\dots,n-1\}$ was
arbitrary, (ii) follows.

Statement (iii) follows trivially from statements (i) and (ii).
\medskip

\noindent {\bf Lemma 2.8.} {\sl Suppose $n\in\N$, and $\mu$ is a
probability measure on $\R^n$ that satisfies Condition $\cH(n)$.
Then the following statements hold:

(A) The probability measure $\theta(\mu)$ on $\R^{Ln}$ (see step
2.4(D)) satisfies Condition $\cH(Ln)$.

(B) Suppose $W$ is an $\R^{Ln}$--valued random vector of the form
$$W:= \left(W_0,W_1,\dots,W_{Ln-1}\right) := \left\langle
\zeta^{(0)}\mid \zeta^{(1)} \mid \zeta^{(2)} \mid \cdots \mid
\zeta^{(L-1)}\right\rangle \eqno{(2.7)}$$
(see (2.4)) where (i)~for each $\ell \in \{0,1,\dots,L-1\}$,
$\zeta^{(\ell)}$ is an $\R^n$--valued random vector such that
$\cL(\zeta^{(\ell)}) = \mu$, and (ii)~these random vectors
$\zeta^{(\ell)}$, $\ell\in \{0,1,\dots,L-1\}$ are independent of each
other. Then (see (2.5)) for every $j\in \{0,1,\dots,L-1\}$,
$\cL(\psi_{n,j}(W)) = \theta(\mu)$.

(C) Suppose the $\R^{Ln}$--valued random vector $W$ is as in
statement (B) (satisfying all conditions there), and $Y:=
(Y_0,Y_1,\dots,Y_{Ln-1})$ is an $\R^{Ln}$--valued random vector such
that $\cL(Y) = \theta(\mu)$. Suppose $S$ is a nonempty subset of
$\{0,1,\dots,Ln-1\}$. Then the following statements (i), (ii), (iii)
hold:

(i) For each integer $M\in \{1,2,\dots,L-1\}$, $E(\sum_{k\in S}Y_k)^M
= E(\sum_{k\in S} W_k)^M$. (In the case where $M$ is odd, both sides
of that equality are $0$.)

(ii) $E(\sum_{k\in S}Y_k)^L \le E(\sum_{k\in S}W_k)^L$.

(iii) For the case where $S = \{0,1,\dots,Ln-1\}$ itself, one has
(see (2.2))}
$$E(\hbox{sum }\, W)^L - E(\hbox{sum}\,\, Y)^L \ge 2^{-L}\cdot L! \cdot
n^{L/2}. \eqno{(2.8)}$$
\medskip

Statements (A) and (B) and their proofs are taken or adapted from
[9]. Their proofs will be included here because of extra
complications in our context.
\medskip

\noindent{\bf Proof.} The proofs of statements (B), (A), and (C)
will be given in that order.
\medskip

\noindent{\bf Proof of statement (B).} Before the random vector $W$
in statement (B) is brought into the picture, some preliminary work
is needed.

Let $j\in \{0,1,\dots,L-1\}$ be arbitrary but fixed. (See the role of
$j$ in statement (B).)

Let $V'_0, V'_1,\dots,V'_{L-1}$ be independent, identically
distributed $\{-1,1\}$--valued random variables such that $P(V'_0
=-1) =P(V'_0 =1) = 1/2$. Define the random vector $V^* := (V^*_0,
V^*_1,\dots,V^*_{L-1})$ as follows: For each $\ell \in
\{0,1,\dots,L-1\} -\{j\}$, $V^*_\ell := V'_\ell$; and
$$ V^*_j := - \prod_{\ell\in \{0,1,\dots,L-1\}-\{j\}} V^*_\ell.
\eqno{(2.9)}$$

Now the random variables $V^*_\ell$, $\ell \in
\{0,1,\dots,L-1\}-\{j\}$ are independent and identically
distributed, $P(V^*_\ell =-1) =P(V^*_\ell =1) = 1/2$ for each $\ell
\in \{0,1,\dots,L-1\}-\{j\}$, and (2.9) holds. These conditions
together uniquely determine $\cL(V^*)$. Also, the random vector $V$
(with $\cL(V) =\nu$) in Remark~2.5 satisfies the analogs of those
conditions (see Remark 2.5(i)(ii)(v)). It follows that
$$\cL(V^*) = \nu. \eqno{(2.10)}$$

Next, let $\xi^{(0)},\xi^{(1)},\dots,\xi^{(L-1)}$ be independent
$\R^n$--valued random vectors such that (i)~$\cL(\xi^{(\ell)})=\mu$
for each $\ell\in \{0,1,\dots,L-1\}$ and
(ii)~$\sigma(\xi^{(0)},\xi^{(1)},\dots,\xi^{(L-1)})$ is independent
of $\sigma(V')$ (and hence independent of $\sigma(V^*))$ --- see
step 2.2(h). Define the $\R^{Ln}$--valued random vectors $W'$ and
$W^*$ by (see (2.3) and (2.4))
$$\eqalignno{
W' &:= \left(W'_0,W'_1,\dots,W'_{Ln-1}\right) \cr &:= \left\langle
V'_0\varphi_n(\xi^{(0)})\mid V'_1\varphi_n(\xi^{(1)})\mid \cdots
\mid V'_{L-1}\varphi_n(\xi^{(L-1)})\right\rangle &(2.11)\cr}$$ and
$$\eqalignno{
W^* &:= \left(W^*_0,W^*_1,\dots,W^*_{Ln-1}\right) \cr &:=
\left\langle V^*_0 \varphi_n(\xi^{(0)})\mid V^*_1
\varphi_n(\xi^{(1)})\mid \cdots  \mid V^*_{L-1}\varphi_n
(\xi^{L-1)})\right\rangle. &(2.12) \cr}$$ By (2.10), (2.12), and
step 2.4(D),
$$\cL(W^*) = \theta(\mu). \eqno{(2.13)}$$

Since (by hypothesis) $\mu$ satisfies Condition $\cH(n)$, one has by
Definition 2.6(a)(e) and Remark~4 in step 2.3(A) that $\forall$
$\ell\in \{0,1,\dots,L-1\}$,
$$\cL(V'_\ell \varphi_n (\xi^{(\ell)})) =
\cL(\xi^{(\ell)}) =\mu. \eqno{(2.14)}$$ Also, the random vectors
$V'_\ell \varphi_n(\xi^{(\ell)})$, $\ell \in \{0,1,\dots,L-1\}$ are
independent (by the properties above). Hence by (2.11) and (2.7) and
the (other) assumptions in statement (B),
$$\cL(W') =\cL(W). \eqno{(2.15)}$$

Our next task is to compare the random vectors $W'$ and $W^*$. For
that purpose, observe that by (2.14) and Definition 2.6(a) (and the
hypothesis of Lemma 2.8),
$$P\left(\prod^{L-1}_{\ell=0} \hbox{sum}\left(V'_\ell
\varphi_n(\xi^{(\ell)})\right) =0\right) =0,
$$
and for all $\omega\in \Omega$ not in that event,
$$V^*_j(\omega)\varphi_n(\xi^{(j)}(\omega)) =
\cases{-V'_j(\omega)\varphi_n(\xi^{(j)}(\omega)) &if
$\prod^{L-1}_{\ell=0}
\hbox{sum}(V'_\ell(\omega)\varphi_n(\xi^{(\ell)}(\omega)))>0$ \cr
V'_j(\omega) \varphi_n(\xi^{(j)}(\omega)) &if $\prod^{L-1}_{\ell=0}
\hbox{sum}(V'_\ell (\omega) \varphi_n(\xi^{(\ell)}(\omega)))<0$.
\cr} \eqno{(2.16)}$$ To verify (2.16), consider first the case where
$\prod^{L-1}_{\ell=0}\hbox{sum}(V'_\ell(\omega)\varphi_n(\xi^{(\ell)}(\omega)))>0$.
This can be rewritten as $\prod^{L-1}_{\ell=0}[V'_\ell(\omega) \cdot
\, \hbox{sum}\, \varphi_n(\xi^{(\ell)}(\omega))]>0$, or
$[\thinspace\prod^{L-1}_{\ell=0} V'_\ell(\omega)] \cdot [\thinspace
\prod^{L-1}_{\ell=0}\, \hbox{sum}\,
\varphi_n(\xi^{(\ell)}(\omega))]>0$. Since
sum~$\varphi_n(\xi^{(\ell)}(\omega))>0$ for each
$\ell\in\{0,1,\dots,L-1\}$, this forces
$\prod^{L-1}_{\ell=0}V'_\ell(\omega)=1$. That in turn implies
$V'_j(\omega) = \prod_{\ell\not= j}V'_\ell(\omega)$, and hence by
(2.9) and its entire sentence,
$$V^*_j(\omega) = -\prod_{\ell\not= j} V^*_\ell(\omega) =
-\prod_{\ell\not= j} V'_\ell(\omega) = -V'_j(\omega),$$ and hence
(2.16) holds. In the other case, where $\prod^{L-1}_{\ell=0}$
sum$(V'_\ell(\omega)\varphi_n(\xi^{(\ell)}(\omega)))<0$, eq.\ (2.16)
holds by a similar argument.

Now of course $V^*_\ell \varphi_n(\xi^{(\ell)})
=V'_\ell\varphi_n(\xi^{(\ell)})$ for each $\ell\in \{0,1,\dots,L-1\}
- \{j\}$. Hence by (2.5), (2.11), (2.12), and (2.16) and its entire
sentence, $$\psi_{n,j}(W') = W^* \,\,\hbox{a.s.} \eqno{(2.17)}$$

Now by (2.15), (2.17), and (2.13), $\cL(\psi_{n,j}(W)) =
\cL(\psi_{n,j}(W')) = \cL(W^*) = \theta(\mu)$. Since
$j\in\{0,1,\dots,L-1\}$ was arbitrary, statement(B) in Lemma 2.8
holds.
\medskip

\noindent{\bf Proof of statement (A).} This will be a continuation
of the argument above for statement (B). (The integer $j$ in that
argument will not play a role here; as a formality, one can take,
say, $j=0$.) Using (2.12) and (2.13), we shall verify, one by one,
the analogs of conditions (a)--(g) in Definition~2.6 for the
probability measure $\theta(\mu)$ on~$\R^{Ln}$.

By (2.14) and the paragraph after (2.10),
$\cL(\xi^{(0)},\xi^{(1)},\dots,\xi^{(L-1)}) = \mu\times
\mu\times\cdots \times \mu$ (the $L$--fold product measure on
$(\R^n)^L)$, which (by Definition 2.6(a) for $\mu$ and a standard
measure--theoretic argument) is absolutely continuous with respect
to Lebesgue measure on $(\R^n)^L$. Hence by (2.3), Definition 2.6(e)
for $\mu$,  and an elementary argument, $\cL(\varphi_n(\xi^{(0)})$,
$\varphi_n(\xi^{(1)}),\dots,$ $\varphi_n(\xi^{(L-1)}))$ is
absolutely continuous. Hence so is $\cL(W^*)$, by (2.12) and an
elementary argument. Hence by (2.13), the analog of condition (a) in
Definition~2.6 holds for~$\theta(\mu)$.

Next, by (2.14), (2.10), Remark 2.5(i)(ii), and the paragraph after
(2.10), one has that
$$\forall\,\,\ell\in\{0,1,\dots,L-1\},\,\quad\cL(V^*_\ell\varphi_n(\xi^{(\ell)})) = \cL(V'_\ell
\varphi_n(\xi^{(\ell)})) =\mu, \eqno{(2.18)}$$ and also that every
$L-1$ of the random vectors $V^*_\ell\varphi_n(\xi^{(\ell)})$,
$\ell\in \{0,1,\dots,L-1\}$ are independent. Hence by (2.12),
Definition 2.6(b)(c) for $\mu$, and a simple argument, one has that
$\cL(W^*_k) = \lambda_{\unp3}$ for each $k\in \{0,1,\dots,Ln-1\}$,
and that the random vector $W^*$ satisfies $(L-1)$--tuplewise
independence. Hence by (2.13), the analogs of conditions (b) and (c)
in Definition~2.6 hold for~$\theta(\mu)$.

Next, by (2.7) and its entire sentence, together with Definition
2.6(d) for $\mu$, the random variables $|W_i|$, $i\in
\{0,1,\dots,Ln-1\}$ are independent. Hence by (2.15), the same holds
for the random variables $|W'_i|$, $0\le i\le Ln-1$. Now by (2.11),
(2.12), and the entire paragraph containing (2.9), $|W^*_i| =
|W'_i|$ for each $i\in \{0,1,\dots,Ln-1\}$; and hence the random
variables $|W^*_i|$, $0\le i\le Ln-1$ are independent. Hence by
(2.12) and (2.13), the analog of condition (d) in Definition 2.6
holds for~$\theta(\mu)$.

Next, by (2.10), Remark 2.5(iii), and the paragraph after (2.10),
$$\eqalignno{\cL&(V^*, \xi^{(0)},\xi^{(1)},\dots,\xi^{(L-1)}) = \nu\times
\mu\times \mu\times\dots\times \mu \cr
&=\cL(-V^*,\xi^{(0)},\xi^{(1)},\dots,\xi^{(L-1)}). &(2.19) \cr}$$
Also, by (2.12), $W^* =
f(V^*,\xi^{(0)},\xi^{(1)},\dots,\xi^{(L-1)})$ for an obvious Borel
function $f:\R^L \times (\R^n)^L\to\R^{nL}$, and $-W^* =
f(-V^*,\xi^{(0)},\xi^{(1)},\dots,\xi^{(L-1)})$. Hence $\cL(-W^*) =
\cL(W^*)$ by (2.19). Hence by (2.13), the analog of condition (e) in
Definition 2.6 holds for~$\theta(\mu)$.
\medskip

{\it Proof of Definition 2.6(f) for $\theta(\mu)$.} Refer to the
notations in step 2.2(c), and refer to Remark 2 in step 2.3(A).

Let $J\in\{0,1,\dots,Ln-1\}$ be arbitrary but fixed.  Refer to (2.12)
and (2.13). To prove the analog of property (f) in Definition 2.6 for
$\theta(\mu)$, it suffices to show that there exists a permutation
$\tau$ of the set $\{0,1,\dots,Ln-1\}$ such that
$$\tau(0) =J \,\,\hbox{and}\,\, \cL(W^*_\tau) = \cL(W^*).
\eqno{(2.20)}$$

Let $\ell'\in \{0,1,\dots,L-1\}$ and $j'\in \{0,1,\dots,n-1\}$ be
such that
$$J = \ell'n +j'. \eqno{(2.21)}$$
Let $\alpha$ be a permutation of the set $\{0,1,\dots,L-1\}$ such
that $\alpha(0) = \ell'$. Then by (2.19), (2.10), and Remark 2.5(iv),
$$\cL\left(V^*_\alpha,\xi^{(\alpha(0))},\xi^{(\alpha(1))},\dots,\xi^{(\alpha(L-1))}\right)
=\nu\times \mu\times \mu\times \dots\times\mu. \eqno{(2.22)}$$
Applying Definition 2.6(f) for $\mu$ (see (2.14)), let $\beta$ be a
permutation of the set $\{0,1,\dots,n-1\}$ such that $\beta(0) =j'$
and (say) $\cL(\xi^{(\alpha(0))}_\beta) =\mu$. Then by (2.22) and a
trivial argument,
$$\cL\left(V^*_\alpha,\xi^{(\alpha(0))}_\beta,\xi^{(\alpha(1))},
\xi^{(\alpha(2))}, \dots,\xi^{(\alpha(L-1))}\right)
=\nu\times \mu\times \mu\times \dots\times\mu. \eqno{(2.23)}$$

Define the $\R^{Ln}$--valued random vectors $W^{**}$ and $W^{***}$ by
$$\eqalignno{W^{**}&:=
\left(W^{**}_0,W^{**}_1,\dots,W^{**}_{Ln-1}\right) \cr &:=
\left\langle (V^*_\alpha)^{\phantom
*}_0 \varphi_n(\xi^{(\alpha(0))})\mid
(V^*_\alpha)^{\phantom *}_1\varphi_n(\xi^{(\alpha(1))})\mid \cdots
\mid (V^*_\alpha)^{\phantom
*}_{L-1}\varphi_n(\xi^{(\alpha(L-1))})\right\rangle &(2.24) \cr}$$ and
$$\eqalignno{W^{***}&:=
\left(W^{***}_0,W^{***}_1,\dots,W^{***}_{Ln-1}\right) \cr &:=
\Bigl\langle (V^*_\alpha)^{\phantom
*}_0\varphi_n(\xi^{(\alpha(0))}_\beta)\mid
(V^*_\alpha)^{\phantom *}_1\varphi_n(\xi^{(\alpha(1))})\mid
(V^*_\alpha)^{\phantom *}_2\varphi_n(\xi^{(\alpha(2))})\mid \cr
&\qquad\qquad\qquad\qquad \qquad\qquad\qquad\qquad\cdots \mid
(V^*_\alpha)^{\phantom
*}_{L-1}\varphi_n(\xi^{(\alpha(L-1))})\Bigr\rangle. &(2.25)
\cr}$$ By (2.19) (its first equality), (2.22), (2.23), (2.24),
(2.25), (2.12), and (2.13),
$$\cL(W^{***}) = \cL(W^{**}) = \cL(W^*) = \theta(\mu).
\eqno{(2.26)}$$

By (2.21) and the sentence after it, together with (2.24) and
(2.12), the random variables $W^{**}_k$, $0\le k\le Ln-1$ in (2.24)
are (to put this informally) simply a permutation of the random
variables $W^*_k$, $0\le k \le Ln-1$ in (2.12), such that
$W^{**}_{j'} = W^*_J$. Similarly, by (2.24), (2.25), Remark 2 in
step 2.3(A), and the sentence after (2.22), the random variables
$W^{***}_k$, $0\le k\le Ln-1$ in (2.25) are simply a permutation of
the random variables $W^{**}_k$, $0 \le k\le Ln-1$ in (2.24), such
that $W^{***}_0 = W^{**}_{j'}$. Hence the random variables
$W^{***}_k$, $0\le k \le Ln-1$ in (2.25) are simply a permutation of
the random variables $W^*_k$, $0\le k\le Ln-1$ in (2.12), such that
$W^{***}_0 = W^*_J$. Thus for the resulting permutation $\tau$ of
the indices $\{0,1,\dots,Ln-1\}$, one has that $W^{***} = W^*_\tau$
and by (2.26), equation (2.20) holds. That completes the argument
for Definition 2.6(f) for $\theta(\mu)$.
\medskip

{\it Proof of Definition 2.6(g) for $\theta(\mu)$.} Let $S$ be an
arbitrary fixed subset of $\{0,1,\dots,Ln-1\}$ such that card~$S=L$.
Refer to (2.12) and (2.13). Our task is to show that $E(\prod_{k\in
S} W^*_k) \le 0$.

For each $\ell\in \{0,1,\dots,L-1\}$, define the set
$$\Phi(\ell):= \{\ell n,\ell n+1,\dots, \ell n +n-1\}.
\eqno{(2.27)}$$
These sets form a partition of the set
$\{0,1,\dots,Ln-1\}$. The argument here will be divided into three
cases according to how many indices $\ell\in \{0,1,\dots,L-1\}$ are
such that the set $S\cap \Phi(\ell)$ is nonempty. First recall from
(2.12) that for each $\ell\in \{0,1,\dots,L-1\}$,
$$\left(W^*_{\ell n},W^*_{\ell n+1},\dots,W^*_{\ell n+n-1}\right) =
V^*_\ell \varphi_n(\xi^{(\ell)}). \eqno{(2.28)}$$

{\it Case 1.} $S\subset \Phi(\ell)$ for some $\ell \in
\{0,1,\dots,L-1\}$. Then by (2.28) and (2.27) (for that $\ell$),
(2.18), and Definition 2.6(g) for $\mu$, $E(\prod_{k\in S} W^*_k) \le
0$.
\medskip

{\it Case 2.} The set $T:= \{\ell \in \{0,1,\dots,L-1\}:S \cap
\Phi(\ell) \not= \emptyset\}$ satisfies $2\le\card\,\, T\le L-1$.
This forces card$(S\cap \Phi(\ell))\le L-1$ for each $\ell \in T$.
For each $\ell\in T$, by (2.28), (2.27), (2.18), and Definition
2.6(b)(c) for $\mu$, $E(\prod_{k\in S\cap \Phi(\ell)}W^*_k) =
\prod_{k\in S\cap \phi(\ell)}EW^*_k =0$. Also, by (2.28), (2.27),
and the phrase immediately after (2.18), the $\sigma$--fields
$\sigma(W^*_k$, $k\in S\cap \Phi(\ell))$, $\ell\in T$ are
independent. It follows that $E(\prod_{k\in S}W^*_k) =0$.
\medskip

{\it Case 3.} For each $\ell\in \{0,1,\dots,L-1\}$, the set
$\Phi(\ell)$ contains exactly one element of $S$. (Recall that card
$S=L$.) For each $\ell\in \{0,1,\dots,L-1\}$, let $k(\ell)$ denote
the element of $S\cap \Phi(\ell)$. By (2.12), for each $\ell\in
\{0,1,\dots,L-1\}$, $W^*_{k(\ell)} = V^*_\ell T_\ell$ where $T_\ell$
is one of the coordinates of the random vector
$\varphi_n(\xi^{(\ell)})$. For each $\ell \in \{0,1,\dots,L-1\}$,
$ET_\ell >0$ by (2.14) and Lemma 2.7(i)(iii). Hence by (2.10), the
paragraph after (2.10), and Remark 2.5(v),
$$\eqalign{
E\left(\prod_{k\in S}W^*_k\right) &= E\left( \prod^{L-1}_{\ell=0}
W^*_{k(\ell)}\right) = E\left(\prod^{L-1}_{\ell=0} V^*_\ell\right)
\cdot \prod^{L-1}_{\ell=0} ET_\ell \cr &= -1\cdot
\prod^{L-1}_{\ell=0}ET_\ell<0. \cr}$$ That completes the argument
for Case~3, for Definition 2.6(g) for $\theta(\mu)$, and for
statement~(A).
\medskip

\noindent{\bf Proof of statement (C).} In this argument, for any
$\R^{Ln}$--valued random vector $\eta:=
(\eta_0,\eta_1,\dots,\eta_{Ln-1})$ and any nonempty set $Q\subset \{
0,1,\dots,Ln-1\}$, the notation $(\eta_k$, $k\in Q)$ will sometimes
be used to denote the random vector
$(\eta_{k(1)},\eta_{k(2)},\dots,\eta_{k(m)})$ where $m= \card\, Q$
and $k(1),k(2),\dots,k(m)$ are in strictly increasing order the
elements of~$Q$.

We shall use appropriate arguments from the proofs of (A) and (B).

For any nonempty set $Q\subset \{0,1,\dots,Ln-1\}$ such that card
$Q\le L-1$, one has (see step 2.4(A)) that
$$\cL(W_k,\, k\in Q) = \lambda^{[m]}_{\unp3} = \cL(Y_k, \, k\in Q).
\eqno{(2.29)}$$ The first equality holds by the entire sentence
containing (2.7), together with Definition 2.6(b)(c) for $\mu$ and a
simple argument. The second equality holds by Definition 2.6(b)(c)
for $\theta(\mu)$ (recall statement (A), proved above).

Refer to the set $S$ in statement (C). To prove item (i) in
statement (C), note that if $1\le M\le L-1$, then for any $M$--tuple
$(k(1),k(2),\dots,K(M))\in S^M$ (the $k(i)$'s need not be distinct),
$E(\prod^M_{i=1}W_{k(i)}) = (\prod^M_{i=1}Y_{k(i)})$ by (2.29) and
its entire sentence. Adding up over all $M$--tuples $\in S^M$, one
obtains $E(\sum_{k\in S}W_k)^M = E(\sum_{k\in S}Y_k)^M$, the
equality in (i). The last sentence in (i) follows from Definition
2.6(e) for~$\theta(\mu)$ (applied to $Y$).
\medskip

{\it Proof of statement (C)(ii).} Let $(k(1),k(2),\dots,k(L))$ be an
arbitrary fixed $L$--tuple $\in S^L$. (The $k(i)$'s need not be
distinct.) In order to prove statement (C)(ii), it suffices to prove
that
$$E\left(\prod^L_{i=1}Y_{k(i)}\right) \le
E\left(\prod^L_{i=1}W_{k(i)}\right). \eqno{(2.30)}$$
(For then by adding up both sides over all $L$--tuples $\in S^L$, one
would obtain statement (C)(ii).)

If two or more of the integers $k(i)$ are equal, then (2.30) holds
with equality by (2.29) and its entire sentence. Therefore, we now
assume that the integers $k(i)$ are distinct.

Now refer to the sets $\Phi(\ell)$, $\ell\in \{0,1,\dots,L-1\}$ in
(2.27). Suppose first that $\{k(1),k(2)$, $\dots,k(L)\}\subset
\Phi(\ell')$ for some $\ell'\in \{0,1,\dots,L-1\}$. By (2.27) and
(2.7) and its entire sentence,
$$\cL\left(W_k,\, k\in \Phi(\ell')\right) =
\cL\left(\zeta^{(\ell')}\right) =\mu. \eqno{(2.31)}$$
Also, since
$\cL(Y) = \theta(\mu) = \cL(W^*)$ by (2.13) and the assumption on
$Y$, one has by (2.27), (2.12), and (2.18),
$$
\cL(Y_k,\, k\in \Phi(\ell')) = \cL(W^*_k,\, k\in \Phi(\ell'))
=\cL(V^*_{\ell'}\varphi_n(\zeta^{(\ell')})) = \mu. \eqno{(2.32)} $$
Hence $\cL(Y_{k(1)}, Y_{k(2)},\dots,Y_{k(L)}) =
\cL(W_{k(1)},W_{k(2)},\dots,W_{k(L)})$, and (2.30) holds with
equality.

Now consider the remaining case where $\{k(1),k(2),\dots,k(L)\}\cap
\Phi(\ell)$ is nonempty for at least two indices $\ell\in
\{0,1,\dots,L-1\}$. Then none of the sets $\Phi(\ell)$, $\ell\in
\{0,1,\dots,L-1\}$ can contain more than $L-1$ of the integers
$k(i)$. Now (2.31) holds for any $\ell'\in \{0,1,\dots,L-1\}$; and
hence by Definition 2.6(b)(c) for $\mu$, one has that for any
$\ell\in \{0,1,\dots,L-1\}$ such that $\{k(1),\dots,k(L)\}\cap
\Phi(\ell)$ is nonempty, $E\prod W_{k(i)} =0$ where the product is
taken over all $i\in \{0,1,\dots,L\}$ such that $k(i) \in
\Phi(\ell)$. Also, by (2.27) and (2.7) and its entire sentence, the
$\sigma$--fields $\sigma(W_k$, $k\in \Phi(\ell))$, $\ell\in
\{0,1,\dots,L-1\}$ are independent. It follows that $E(\prod^L_{i=1}
W_{k(i)}) =0$. Since (by hypothesis) $\cL(Y) = \theta(\mu)$, one has
by Definition 2.6(g) for $\theta(\mu)$ that
$E(\prod^L_{i=1}Y_{k(i)})\le 0$. Thus (2.30) holds. That complets
the proof of statement (C)(ii).
\medskip

{\it Proof of statement (C)(iii).} Refer again to (2.27). For each
$\ell\in \{0,1,\dots,L-1\}$, define the random variables
$$\eta'_{\ell} := \sum_{k\in \Phi(\ell)}W_k
\quad\hbox{and}\quad \eta^*_\ell := \sum_{k\in \Phi(\ell)}Y_k.
\eqno{(2.33)}$$ Then
$$\hbox{sum}\,W = \sum^{L-1}_{\ell=0} \eta'_\ell \quad\hbox{and$\quad$ sum
 } Y = \sum^{L-1}_{\ell=0} \eta^*_\ell. \eqno{(2.34)}$$

By (2.27) and (2.7) and its entire sentence (note that $\eta'_\ell
=$ sum $\zeta^{(\ell)}$ for each $\ell\in \{ 0,1,\dots,L-1\})$, the
random variables $\eta'_\ell$, $\ell\in \{0,1,\dots,L-1\}$ are
independent. Recall again that $\cL(W^*) = \theta(\mu) = \cL(Y)$ by
(2.13) and the assumption on $Y$. By (2.33), (2.27), (2.12), and the
phrase immediately after (2.18), every $L-1$ of the random variables
$\eta^*_\ell$, $\ell\in \{0,1,\dots,L-1\}$ are independent. Hence by
(2.31), (2.32), and (2.33), if $j(1),j(2),\dots,j(L-1)$ are distinct
elements of $\{0,1,\dots,\L-1\}$, then
$$\cL\left(\eta'_{j(1)},\eta'_{j(2)},\dots,\eta'_{j(L-1)}\right) =
\cL\left( \eta^*_{j(1)},\eta^*_{j(2)},\dots,\eta^*_{j(L-1)}\right).
$$
Hence, if $j(1),j(2), \dots,j(L)$ are each an element of
$\{0,1,\dots,L-1\}$ and two or more of the $j(i)$'s are equal, then
$E(\prod^L_{i=1}\eta'_{j(i)}) = E(\prod^L_{i=1}\eta^*_{j(i)})$.
Hence by (2.34) and a simple calculation,
$$E(\hbox{sum}\, W)^L - E(\hbox{sum}\, Y)^L = L!\cdot
\left[E\left(\prod^{L-1}_{\ell=0} \eta'_\ell\right) -
E\left(\prod^{L-1}_{\ell=0} \eta^*_\ell\right)\right].
\eqno{(2.35)}$$

Now by (2.7) and its entire sentence, and Definition 2.6(b) for
$\mu$, one has that $EW_k =0$ for each $k\in \{0,1\dots, Ln-1\}$; and
hence by (2.33), $E\eta'_\ell =0$ for each $\ell\in
\{0,1,\dots,L-1\}$. Hence by the sentence after (2.34),
$E(\prod^{L-1}_{\ell=0} \eta'_\ell) =0$. Hence by (2.35),
$$E(\hbox{sum}\, W)^L -E(\hbox{sum }\, Y)^L = -L! \cdot
E\left(\sum^{L-1}_{\ell=0} \eta^*_\ell\right). \eqno{(2.36)}$$

Now recall again the equality $\cL(W^*) =\cL(Y)$ in the second
sentence after (2.34). By (2.33), (2.27), and (2.12),
$$
E\left(\prod^{L-1}_{\ell=0} \eta^*_\ell\right) = E\left[
\prod^{L-1}_{\ell=0} \hbox{sum}\left(V^*_\ell
\varphi_n(\xi^{(\ell)})\right)\right] = E\left[\prod^{L-1}_{\ell=0}
\left(V^*_\ell \cdot
\hbox{sum}\left(\varphi_n(\xi^{(\ell)})\right)\right)\right]. $$
Hence by (2.10), Remark 2.5(v), and the sentence after (2.10),
$$\eqalignno{
E\left( \prod^{L-1}_{\ell=0} \eta^*_\ell\right) &=
E\left[\prod^{L-1}_{\ell=0} V^*_\ell\right] \cdot
\prod^{L-1}_{\ell=0} E\left[ \hbox{sum}\, \varphi_n
(\xi^{(\ell)})\right] \cr &= -1 \cdot \left[ E\left(\hbox{sum}\,
\varphi_n (\xi^{(0)})\right)\right]^L. &(2.37) \cr}$$ Now by the
sentence after (2.10), together with Lemma 2.7(i), $E(\hbox{sum}\,
\varphi_n(\xi^{(0)}))\ge (1/2)n^{1/2}$. Hence by (2.37),
$E(\prod^{L-1}_{\ell=0} \eta^*_\ell) \le -2^{-L}n^{L/2}$. Hence by
(2.36), eq.\  (2.8) holds. That completes the proof of statement
C(iii), and of Lemma 2.8.
\medskip

\noindent{\bf Lemma 2.9.} {\sl Suppose $n$ is a positive integer,
and $\mu$ is a probability measure on $\R^n$ that satisfies
Condition $\cH(n)$ (see Definition 2.6). Suppose $W:= (W_k$,
$k\in\Z)$ and $Y:= (Y_k$, $k\in\Z)$ are random sequences such that
(see step 2.4(B)(D)) $\cL(W) = \cD(n,\mu)$ and $\cL(Y) =
\cD(Ln,\theta(\mu))$. Then the following statements hold:

(A) The random variables $W_k$, $k\in\Z$ satisfy $(L-1)$--tuplewise
independence as well as (for all $k\in\Z$) $\cL(W_k) =
\lambda_{\unp3}$ (see step 2.4(A)); and the same holds for the
random variables $Y_k$, $k\in\Z$.

(B) For any nonempty finite set $S\subset \Z$, one has that
$E(\sum_{k\in S} Y_k)^L \le E(\sum_{k\in S}W_k)^L$.

(C) If $j\in \Z$, and $S$ is a finite subset of $\Z$ such that
$\{j,j+1,j+2,\dots,j+2Ln-1\}\subset S$, then}
$$E\left( \sum_{k\in S}W_k\right)^L - E\left(\sum_{k\in
S}Y_k\right)^L \ge 2^{-L}\cdot L! \cdot n^{L/2}. \eqno{(2.38)}$$
\medskip

\noindent{\bf Proof.} Statement (A) holds by Lemma 2.8(A), Definition
2.6(b)(c) (for both $\mu$ and $\theta(\mu))$, step 2.4(B), and a
trivial extra argument.

Statements (B) and (C) will be proved together. For that purpose,
define for each $u\in\Z$ the set
$$\Lambda(u) := \{Lnu, Lnu+1,Lnu+2,\dots,Lnu+Ln-1\}. \eqno{(2.39)}$$
In this argument, we shall follow the convention in the first
paragraph of the proof of statement (C) in Lemma 2.8.  By step
2.4(B) and a trivial argument, the random vectors $(W_k$, $k\in
\Lambda(u))$, $u\in\Z$ are independent of each other; and the random
vectors $(Y_k$, $k\in \Lambda(u))$, $u\in\Z$ are independent of each
other.

Now suppose $S$ is a nonempty finite subset of $\Z$. Let $Q$ denote
the (nonempty, finite) set of all integers $u$ such that the set
$S\cap \Lambda(u)$ is nonempty. For each $q\in Q$, define the random
variables
$$T'_q := \sum_{k\in S\cap \Lambda(q)} W_k \quad\hbox{and}\quad T^*_q :=
\sum_{k\in S\cap \Lambda(q)}Y_k. \eqno{(2.40)}$$ Then
$$\sum_{k\in S} W_k = \sum_{q\in Q}T'_q \quad\hbox{and}\quad \sum_{k\in
S}Y_k = \sum_{q\in Q} T^*_q. \eqno{(2.41)}$$

Now for each $q\in \Z$, by step 2.4(B), the random vectors
$(W_{Lnq},W_{Lnq+1},\dots, W_{Lnq+Ln-1})$ and
$(Y_{Lnq},Y_{Lnq+1},\dots,Y_{Lnq+Ln-1})$ satisfy the hypothesis of
Lemma 2.8(C). Hence for each $q\in Q$ and each $m\in
\{1,2,\dots,L-1\}$, by (2.40) and Lemma 2.8(C)(i), $E(T'_q)^m =
E(T^*_q)^m$. Also, the random variables $T'_q$, $q\in Q$ are
independent, and the random variables $T^*_q$, $q\in Q$ are
independent, by (2.40), (2.39), and the two sentences right after
(2.39). Hence, if a given $L$--tuple $(q(1),q(2),\dots,q(L))\in Q^L$
is such that at least two of the entries $q(i)$ are different from
each other, then by a simple argument, $E(\prod^L_{i=1}T'_{q(i)}) =
E(\prod^L_{i=1}T^*_{q(i)})$. Hence by (2.41) and a simple
calculation,
$$E\left[ \sum_{k\in S}W_k\right]^L - E\left[\sum_{k\in S}Y_k\right]^L
= \sum_{q\in Q}\left[ E(T'_q)^L - E(T^*_q)^L\right]. \eqno{(2.42)}$$

Now by (2.40), the sentence after (2.41), and Lemma 2.8(C)(ii),
$E(T'_q)^L \ge E(T^*_q)^L$ for every $q\in Q$. Hence by (2.42),
statement (B) (in Lemma 2.9) holds.

To continue this argument in order to prove statement (C), note that
under the hypothesis of (C), there exists $p\in Q$ such that (see
(2.39) again) $\Lambda(p)\subset S$. By (2.40), the sentence after
(2.41), and Lemma 2.8(C)(iii), $E(T'_p)^L -E(T^*_p)^L \ge 2^{-L}
\cdot L! \cdot n^{L/2}$. Hence by (2.42) and the sentence right
after it, equation (2.38) holds. That completes the proof of
statement (C), and of Lemma 2.9.
\bigskip

\noindent{\bf 3. Part 2 of proof of Theorem 1.1}
\medskip

This section is a direct continuation of section 2.  As in section
2, the argument here in section 3 will be divided into several
``steps,'' including some ``lemmas.'' \medskip

 \noindent{\bf Step
3.1.} Refer again to the even integer $L\ge 6$ fixed in step 2.1.
Refer to steps 2.2(a) and 2.4(A)(D). Recursively define as follows,
for each $n\in \{0,1,2,\dots\,\}$, the probability measure $\mu_n$
on $\R^{L\uparrow n}$: (i)~$\mu_0:= \lambda_{\unp3}$, and (ii)~for
each $n\ge 0$, $\mu_{n+1} := \theta(\mu_n)$.
\medskip

{\it Remark.} For each $n\in \{0,1,2,\dots\,\}$, the probability
measure $\mu_n$ satisfies Condition $\cH(L^n)$. (For $n=0$, that
holds by Definition 2.6 and a trivial argument; and then for $n\ge
1$, it holds by Lemma 2.8(A) and induction.)
\medskip

For the next lemma, and for the rest of section 3 here, let $Z$ be a
$N(0,1)$ random variable.
\medskip

\noindent{\bf Lemma 3.2.} {\sl Suppose $h$ and $m$ are positive
integers such that $m\le h$. Suppose $Y:= (Y_k$, $k\in \Z)$ is a
random sequence such that $\cL(Y) = \cD(L^h,\mu_h)$ (see step
2.4(B)). Suppose $S$ is a finite set $\subset \Z$ such that for some
integer $j$, one has that $\{j,j+1,j+2,\dots,j+2L^m-1\}\subset S$.
Then}
$$E\left[ \sum_{k\in S} Y_k\right]^L \le E\left[ (\card \,
S)^{1/2}Z\right]^L -2^{-L}\cdot L! \cdot L^{(m-1)L/2}. \eqno{(3.1)}$$

\noindent{\bf Proof.} For each $n\in \{0,1,2,\dots\, \}$, let
$W^{(n)}:= (W^{(n)}_k$, $k\in \Z)$ be a random sequence such that
$\cL(W^{(n)}) = \cD(L^n$, $\mu_n)$. Then by Lemma 2.9(B),
$$\forall\, n\ge 1,\quad E\left[ \sum_{k\in S}W^{(n)}_k \right]^L \le
E\left[ \sum_{k\in S}W^{(n-1)}_k\right]^L. \eqno{(3.2)}$$ Also by
Lemma 2.9(C) (with $n=L^{m-1}$) and the assumptions on $S$,
$$E\left[ \sum_{k\in S}W_k^{(m)}\right]^L \le E\left[ \sum_{k\in
S}W^{(m-1)}_k\right]^L -2^{-L} \cdot L! \cdot L^{(m-1)L/2}.
$$
Combining that with (3.2) and trivial induction, one has that
$$\forall\, n\ge m,\quad E\left[ \sum_{k\in S} W^{(n)}_k\right]^L \le
E\left[ \sum_{k\in S}W^{(0)}_k\right]^L -2^{-L}\cdot L!\cdot
L^{(m-1)L/2}.$$ Since $h\ge m$ (by hypothesis) and $\cL(Y)  =
\cL(W^{(h)})$, it follows that
$$E\left[ \sum_{k\in S}Y_k\right]^L \le E\left[ \sum_{k\in
S}W^{(0)}_k\right]^L -2^{-L}\cdot L!\cdot L^{(m-1)L/2}.$$ Hence, to
complete the proof of (3.1), it suffices to prove that
$$E\left[ \sum_{k\in S}W^{(0)}_k\right]^L \le E\left[ (\card \,
S)^{1/2}Z\right]^L. \eqno{(3.3)}$$

Now (see steps 3.1 and 2.4(A) again), the random variables
$W^{(0)}_k$, $k\in \Z$ are independent and identically distributed,
with $\cL(W^{(0)}_0) = \lambda_{\unp3}$. Let $Z_k$, $k\in\Z$ be
independent $N(0,1)$ random variables. If $p$ is a positive integer,
$k(1),k(2),\dots,k(p)$ are distinct  integers, and
$a(1),a(2),\dots,a(p)$ are (not necessarily distinct) positive
integers, then by the Remark in step 2.4(A),
$$E\left[
\prod^p_{i=1}\left(W^{(0)}_{k(i)}\right)^{a(i)}\right] =
\prod^p_{i=1} E\left(W^{(0)}_{k(i)}\right)^{a(i)}  \le \prod^p_{i=1}
EZ^{a(i)}_{k(i)} = E\left[ \prod^p_{i=1}Z^{a(i)}_{k(i)}\right]. $$
Thus for any $L$--tuple $(k(1),k(2),\dots,k(L))\in \Z^L$ (the
$k(\ell)$'s need not be distinct), \break \noindent $E[\thinspace
\prod^L_{\ell=1}W^{(0)}_{k(\ell)}] \le
E[\thinspace\prod^L_{\ell=1}Z_{k(\ell)}]$. Adding up both sides over
all $L$--tuples $\in S^L$, one obtains $E[\thinspace\sum_{k\in
S}W^{(0)}_k]^L \le E[\thinspace\sum_{k\in S}Z_k]^L$. Since
$\cL(\sum_{k\in S}Z_k) = N(0,\card \, S) = \cL((\card\, S)^{1/2}Z)$,
eq.\ (3.3) follows. That completes the proof of (3.1) and of Lemma
3.2.
\medskip

\noindent{\bf Step 3.3.} Let $\kappa_1,\kappa_2,\kappa_3,\dots$ be a
sequence of independent, identically distributed random variables,
taking their values in the set $\{0,1,2,\dots,L-1\}$ and uniformly
distributed on that set.

Define the random variables $J(n)$, $n\in \{0,1,2,\dots\,\}$ as
follows:
$$J(0):=0\, \,\, \hbox{(constant), and for each $n\in\N$,}\,\, J(n) :=
\sum^n_{u=1}L^{u-1}\kappa_u. \eqno{(3.4)}$$
\medskip

{\it Remarks.} Recall that in the base--$L$ number system, for a
given $n\ge 1$ and a given $\beta\in \{0,1,2,\dots,L^{n}-1\}$, there
is a unique representation $\beta = \sum^n_{u=1}L^{u-1}\alpha_u$
such that $\alpha_1,\alpha_2,\dots,\alpha_n\in \{0,1,\dots,L-1\}$.
Also, by (3.4),
$$\forall\, n\in\N,\,\,J(n) = L^{n-1} \kappa_n +J(n-1).
\eqno{(3.5)}$$
From those facts and trivial arguments, one has the
following observations:
\medskip

{\it Remark 1.} For each $n\in \{0,1,2,\dots\,\}$, the random
variable $J(n)$ takes its values in the set $\{0,1,2,\dots,L^n-1\}$
and is uniformly distributed on that set.
\medskip

{\it Remark 2.} For each $n\in \N$ (and each $\omega\in \Omega$),
$J(n) \ge J(n-1)$ and in fact
$$\eqalignno{
-J(n) &\le -J(n-1)\le 0 \cr &\le -J(n-1) +L^{n-1} -1 \le -J(n) +L^n
-1. &(3.6)\cr}$$ (For the last inequality, apply (3.5) and the fact
$L^{n-1}\kappa_n \le L^{n-1}(L-1) = L^n -L^{n-1}$.)
\medskip

{\it Remark 3.} For each $n\in\N$ and each $\omega\in \Omega$ such
that $1\le \kappa_n(\omega) \le L-2$, one has (again by (3.5)) that
$$\eqalignno{
&-J(n)(\omega) \le -L^{n-1} - J(n-1)(\omega) \quad\hbox{and} \cr
&(-J(n-1)(\omega) + L^{n-1}-1) + L^{n-1} \le (-J(n)(\omega) +L^n-1).
&(3.7) \cr}$$
\medskip

{\it Remark 4.} For each $n\in\N$ and each $j\in
\{0,1,2,\dots,L^n-1\}$, there exists a unique choice of integers
$j_0 \in \{0,1,2,\dots,L^{n-1}-1\}$ and $k_0 \in \{0,1,\dots,L-1\}$
such that
$$\{J(n) =j\} = \{J(n-1) = j_0\}\cap \{\kappa_n = k_0\}.
\eqno{(3.8)}$$

{\it Remark 5.} For each $n\in\N$, the two random variables $J(n-1)$
and $\kappa_n$ are independent.
\medskip

\noindent {\bf Step 3.4.} For each $n\in \{0,1,2,\dots\,\}$, we shall
construct a sequence $X^{(n)}:= (X^{(n)}_k$, $k\in \Z)$ of random
variables. The definition will be recursive, and is as follows:
\medskip

(A) To start off, let $X^{(0)}:= (X^{(0)}_k$, $k\in \Z)$ be a
sequence of independent, identically distributed random variables
such that (i)~$\cL(X^{(0)}_0) = \lambda_{\unp3}$ (again see step
2.4(A)) and (ii)~this sequence $X^{(0)}$ is independent of the
sequence $(\kappa_1,\kappa_2,\kappa_3,\dots\,)$ in step 3.3.

(B) Now suppose $n\ge 1$, and the random sequence $X^{(n-1)}:=
(X^{(n-1)}_k$, $k\in \Z)$ is already defined. Define the random
sequence $X^{(n)}:= (X^{(n)}_k$, $k\in \Z)$ as follows:
\medskip

For a given $\omega\in \Omega$, in the notations in step 2.2(a)(d)
and step 2.3(C),
$$\eqalignno{
&(X^{(n)}_k(\omega): -J(n)(\omega) \le k \le -J(n)(\omega) +L^n-1)
\cr &:= \cases{\psi_{L\uparrow (n-1),1}(X^{(n-1)}_k(\omega):
-J(n)(\omega) \le k\le -J(n)(\omega) +L^n-1) & \cr \qquad\qquad
\hbox{if $J(n)(\omega) = J(n-1)(\omega)$} & \cr
\psi_{L\uparrow(n-1),0} (X^{(n-1)}_k(\omega): -J(n)(\omega) \le k\le
-J(n)(\omega) +L^n-1) & \cr \qquad\qquad \hbox{if
$J(n)(\omega)>J(n-1)(\omega)$,} \cr} \cr &&(3.9) \cr}$$ and for
every integer $K\in \Z-\{0\}$,
$$\eqalignno{
&(X^{(n)}_k(\omega): -J(n)(\omega) + KL^n \le k \le -J(n)(\omega)
+KL^n +L^n-1) \cr &:= \psi_{L\uparrow (n-1),0} (X^{(n-1)}_k(\omega):
 -J(n)(\omega) +KL^n \le k\le -J(n)(\omega) +KL^n +L^n-1). \cr
&&(3.10)\cr}$$

That completes the recursive definition. By induction on $n$, one has
(see step 2.2(h)) that
$$\forall\, n\in \N,\,\, \sigma(X^{(n)}) \subset \sigma(X^{(0)}) \vee
\sigma(\kappa_1,\kappa_2,\dots,\kappa_n). \eqno{(3.11)}$$ By step
3.3 and paragraph (A) (property (ii)) above, one now has that for
all $n\in \{0,1,2,\dots\,\}$, the two $\sigma$--fields
$\sigma(J(n)$, $X^{(n)})$ and $\sigma(\kappa_{n+1})$ are
independent. \medskip

 \noindent{\bf Lemma 3.5.} {\sl For each $n\in
\{0,1,2,\dots\,\}$, one has that }
$$\forall\,\, j\in
\{0,1,2,\dots,L^n-1\},\quad\cL\left(X^{(n)}_{-j+k},\, k\in \Z\mid
J(n)=j\right) =\cD\left(L^n,\, \mu_n\right). \eqno{(3.12)}$$

\noindent{\bf Proof.} Eq.\ (3.12) holds for $n=0$ by (3.4) and step
3.4(A), step 3.1, and step 2.4(B).

Now for induction, suppose that $n\ge 1$, and that (3.12) holds with
$n$ replaced by $n-1$. To complete the induction step and the proof,
our task is to prove (3.12) for the given~$n$.

Referring to (3.12) suppose $j\in \{0,1,2,\dots,L^n-1\}$.

Referring to step 3.3 (its Remark 4), let $j_0\in
\{0,1,\dots,L^{n-1}-1\}$ and $k_0 \in \{0,1,\dots,L-1\}$ be such that
(3.8) holds. By the sentence after (3.11), the $\sigma$--fields
$\sigma(J(n-1),X^{(n-1)})$ and $\sigma(\kappa_n)$ are independent.
Hence by (3.8) and a simple calculation,
$$\eqalignno{
\forall\, B\in \sigma(X^{(n-1)}), \quad P(B\mid J(n) =j) &= P(B\mid
\{J(n-1) = j_0\}\cap \{\kappa_n = k_0\}) \cr &= P(B\mid J(n-1) =
j_0). &(3.13) \cr}$$

Now by (3.8) and (3.5), $j= L^{n-1}k_0 +j_0$. By the induction
hypothesis,
$$\cL\left(X^{(n-1)}_{-j(0) +k}, \, k\in \Z \mid J(n-1) = j_0\right)
= \cD\left(L^{n-1},\mu_{n-1}\right) \eqno(3.14) $$ (where $j_0$ is
also written as $j(0)$ for typographical convenience). Also (see
step 2.4(B)), the distribution $\cD(L^{n-1}$, $\mu_{n-1})$ on
$\R^\Z$ is invariant under a shift of the indices by (any multiple
of) $L^{n-1}$. Hence by (3.14),
$$\cL\left(X^{(n-1)}_{-j+k},\, k\in \Z \mid J(n-1) = j_0\right) = \cD
\left(L^{n-1},\mu_{n-1}\right).$$
Hence by (3.13),
$$\cL\left(X^{(n-1)}_{-j+k}, \, k\in \Z\mid J(n)=j\right) = \cD
\left( L^{n-1}, \mu_{n-1}\right). \eqno{(3.15)}$$

Referring to step 2.2(d), consider the following three classes of
random vectors for $K\in\Z$:
$$\eqalign{
\zeta'_K &:= \left(X^{(n-1)}_k :-j +KL^{n-1} \le k\le -j+KL^{n-1}
+L^{n-1}-1\right), \cr \zeta''_K &:= \left(X^{(n-1)}_k :-j +KL^n \le
k\le -j +KL^n +L^n -1\right), \cr \zeta'''_K &:= \left( X^{(n)}_k
:-j +KL^n \le k\le -j +KL^n +L^n-1\right). \cr}$$ Referring to
(2.4), note that for each $K\in\Z$, $\zeta''_{K} = \left\langle
\zeta'_{LK} \mid \zeta'_{LK+1} \mid \zeta'_{LK+2}\mid \cdots \mid
\zeta'_{LK+L-1}\right\rangle$. Also, by (3.9)--(3.10) (and the
inclusion $\{J(n) =j\}\subset \{J(n-1) =j_0\}$ from (3.8)), for each
$K\in \Z$, either (if $K\not= 0$, or $K=0$ and $j>j_0$)
$\zeta'''_K(\omega) = \psi_{L\uparrow (n-1),0}(\zeta''_K(\omega))$
for all $\omega\in \{J(n) =j\}$, or (if $K=0$ and $j=j_0$)
$\zeta'''_K(\omega) = \psi_{L\uparrow (n-1),1}(\zeta''_K(\omega))$
for all $\omega\in \{J(n) =j\}$. Now by (3.15) and step 2.4(B),
conditional on the event $\{J(n)=j\}$, the random vectors
$\zeta'_K$, $K\in\Z$ are independent and have distribution
$\mu_{n-1}$. Hence conditional on $\{J(n) =j\}$, the random vectors
$\zeta''_K$, $K\in\Z$ are independent. Hence by step 3.1 and Lemma
2.8(B), conditional on $\{J(n) =j\}$, the random vectors
$\zeta'''_K$, $K\in\Z$ are independent and have distribution
$\mu_n$. Thus by step 2.4(B), eq.\ (3.12) holds for our given $n$
and $j$. Since $j\in \{0,1,\dots,L^n-1\}$ was arbitrary, that
completes the induction step and the proof.
\medskip

\noindent {\bf Lemma 3.6.} {\sl For each $n\in\N$, the random
sequence $X^{(n)}$ is strictly stationary and satisfies
$\cL(X^{(n)}_0) = \lambda_{\unp3}$ (see step 2.4(A)) and
$(L-1)$--tuplewise independence, and the random variables
$|X^{(n)}_k|$, $k\in\Z$ are independent.}
\medskip

\noindent {\bf Proof.} Let $n\in\N$ be arbitrary but fixed.

To prove strict stationarity, let $\zeta := (\zeta_k$, $k\in\Z)$ be
a random sequence such that (see step 2.4(B)) $\cL(\zeta) =
\cD(L^n$, $\mu_n)$. For each $i\in\Z$, define the probability
measure $\eta_i$ on $\R^\Z$ by $\eta_i:= \cL(\zeta_{i+k}$,
$k\in\Z)$. Then $\eta_0 = \cD(\cL^n$, $\mu_n)$, and from step 2.4(B)
one has the ``periodicity'' $\eta_{i+L\uparrow n} = \eta_i$ for
$i\in\Z$. By Lemma 3.5, if $j\in \{0,1,\dots,L^n-1\}$ then
$\cL(X^{(n)}_{-j+k}$, $k\in\Z \mid J(n) =j) =\eta_0$, and if also
$i\in\Z$ then $\cL(X^{(n)}_{-j+k+i}$, $k\in\Z\mid J(n) =j) =\eta_i$.
Hence by step 3.3 (its Remark 1) and a simple calculation, for any
$i\in\Z$,
$$\eqalign{
\cL\left(X^{(n)}_{i+k}, k\in\Z\right)  &= L^{-n} \cdot
\sum^{(L\uparrow n)-1}_{j=0} \cL\left(X^{(n)}_{i+k}, k\in\Z\mid J(n)
=j\right) \cr &= L^{-n} \cdot \sum^{(L\uparrow n)-1}_{j=0}
\cL\left(X^{(n)}_{-j+k+i+j}, k\in\Z\mid J(n) =j\right) \cr &= L^{-n}
\cdot \sum^{(L\uparrow n)-1}_{j=0} \eta_{i+j} = L^{-n} \cdot
\sum^{(L\uparrow n)-1}_{u=0} \eta_{u} \cr}$$ where the last equality
comes from the ``periodicity'' mentioned above. The last expression
does not depend on $i$. Hence the sequence $X^{(n)}$ is strictly
stationary.

Next, by Lemma 3.5, step 3.1 (its Remark), Definition 2.6(b)(c)(d),
and step 2.4(B) (and where necessary a trivial extra argument), one
has the following observations: First,
$$\eqalign{
\cL\left(X^{(n)}_0\right) &= L^{-n} \cdot \sum^{(L\uparrow
n)-1}_{j=0} \cL\left( X^{(n)}_0 \mid J(n) =j\right) \cr &= L^{-n}
\cdot \sum^{(L\uparrow n)-1}_{j=0} \lambda_{\unp3} =
\lambda_{\unp3}. \cr}$$
Second, if $k(1),k(2),\dots,k(L-1)$ are any
choice of $L-1$ distinct integers, then (with a notation at the end
of step 2.4(A))
$$\eqalign{
\cL\left(X^{(n)}_{k(1)},X^{(n)}_{k(2)},\dots,X^{(n)}_{k(L-1)}\right)
\ &= L^{-n} \cdot \sum^{(L\uparrow n)-1}_{j=0} \cL
\left(X^{(n)}_{k(1)},X^{(n)}_{k(2)},\dots,X^{(n)}_{k(L-1)} \mid
J(n)=j\right) \cr &= L^{-n} \cdot \sum^{(L\uparrow n)-1}_{j=0}
\lambda^{[L-1]}_{\unp3}\,\,\, =\,\, \lambda^{[L-1]}_{\unp3}\, ;
\cr}$$ and hence the sequence $X^{(n)}$ satisfies $(L-1)$--tuplewise
independence. Third, by a similar argument, the random variables
$|X^{(n)}_k|$, $k\in\Z$ are independent. The proof of Lemma 3.6 is
complete.
\medskip

\noindent{\bf Lemma 3.7.} {\it Suppose $h$ and $m$ are positive
integers such that $m\le h$. Suppose $p$ is a positive integer such
that $p\ge 2L^m$. Then (recall the sentence preceding Lemma 3.2)}
$$E\left[ \sum^p_{k=1}X^{(h)}_k\right]^L \le E\left[
p^{1/2}Z\right]^L - 2^{-L} \cdot L! \cdot L^{(m-1)L/2}.
\eqno{(3.16)}$$

\noindent{\bf Proof.}  By Lemma 3.5 and Lemma 3.2, for each $j\in
\{0,1,\dots,L^h-1\}$,
$$E\left(\left[ \sum^p_{k=1} X^{(h)}_k\right]^L \,\bigg| \,\,J(n) =j\right)
\le [\hbox{RHS of (3.16)}].$$
Hence by step 3.3 (its Remark 1) and a trivial calculation, (3.16)
holds. \medskip

\noindent{\bf Lemma 3.8.} {\sl There exists an event $G$ with $P(G)
=1$ such that the following holds: For every $\omega\in G$ and every
$k\in \Z$, there exists a positive integer $m = m(\omega,k)$ such
that}
$$\forall\,\, n\ge m, \,\,\, X^{(n)}_k(\omega) =X^{(m)}_k(\omega).
\eqno{(3.17)}$$

\noindent{\bf Proof.} First some preliminary calculations will be
useful.

Suppose $\omega\in \Omega$, $n\ge 1$, and $J(n)(\omega) =
J(n-1)(\omega)$. Then by (3.9) and step 2.3(C), $X^{(n)}_k(\omega) =
X^{(n-1)}_k (\omega)$ will hold for all $k\in \Z$ such that
$-J(n)(\omega) \le k\le -J(n)(\omega) +L^n-1$ with the possible
exception of the indices $k$ such that $-J(n)(\omega) +L^{n-1} \le k
\le -J(n)(\omega) +2L^{n-1}-1$. Thus (if $J(n)(\omega) =
J(n-1)(\omega))$ $X^{(n)}_k(\omega) = X^{(n-1)}_k(\omega)$ will hold
for all $k$ such that $-J(n)(\omega) \le k\le -J(n)(\omega)
+L^{n-1}-1$, that is, for all $k$ such that
$$-J(n-1)(\omega) \le k\le -J(n-1)(\omega) +L^{n-1}-1.
\eqno{(3.18)}$$

Now suppose $\omega\in \Omega$, $n\ge 1$, and (instead)
$J(n)(\omega)>J(n-1)(\omega)$. Then by (3.5), $J(n)(\omega) \ge
J(n-1)(\omega) +L^{n-1}$ and hence $-J(n)(\omega) +L^{n-1} \le
-J(n-1)(\omega)$; and also by (3.9) and step 2.3(C),
$X^{(n)}_k(\omega) = X^{(n-1)}_k(\omega)$ will hold for (at least)
all integers $k$ such that $-J(n)(\omega) +L^{n-1} \le k\le
-J(n)(\omega) +L^n-1$, and in particular (see also the last
inequality in (3.6)) for all integers $k$ such that (3.18) holds.

By the preceding two paragraphs, one has that if $\omega\in\Omega$
and $n\ge 1$, then $X^{(n)}_k(\omega) = X^{(n-1)}_k(\omega)$ for all
integers $k$ such that (3.18) holds.

Hence by (3.6) and induction, one has that if $\omega\in\Omega$,
$n\ge 1$, and $k$ satisfies (3.18), then $X^{(n-1)}_k(\omega) =
X^{(n)}_k(\omega) = X^{(n+1)}_k(\omega) = X^{(n+2)}_k(\omega) =
\dots\,\,$.

Now define the event $G$ by
$$G:= \{ 1\le \kappa_n \le L-2 \,\,\hbox{for infinitely many $n\in \N$}\}.$$
By (2.1) and step 3.3, $P(1\le \kappa_n \le L-2) \ge 2/3$ for each
$n\in\N$, and hence by step 3.3 and the (second) Borel--Cantelli
Lemma, $P(G)=1$. Also, for each $\omega\in G$, [LHS of (3.18)]
$\longrightarrow -\infty$ and [RHS of (3.18)] $\longrightarrow
\infty$ as $n\to \infty$ by (3.7). Hence for every $\omega\in G$ and
every $k\in \Z$, there exists $n\in\N$ such that (3.18) holds. Hence
by the observation in the preceding paragraph, Lemma 3.8 holds.
\medskip

\noindent{\bf Step 3.9.} Define the sequence $X:= (X_k$, $k\in\Z)$ of
random variables as follows: For each $k\in\Z$,
$$X_k := \lim_{n\to\infty} X^{(n)}_k. \eqno{(3.19)}$$
By Lemma 3.8, for each $k\in\Z$, this limit exists almost surely.
\medskip

\noindent{\bf Step 3.10.} The task now is to show that the sequence
$X$ defined in (3.19) has all of the properties stated in
Theorem~1.1. \medskip

Suppose $p$ is a positive integer and $k(1),k(2),\dots,k(p)$ are
distinct elements of $\Z$. By (3.19), for any $p$--tuple
$(t_1,t_2,\dots,t_p) \in\R^p$, one has that $\sum^p_{i=1} t_i
X^{(n)}_{k(i)} \longrightarrow \sum^p_{i=1}t_iX_{k(i)}$ a.s.\ as
$n\to \infty$. Hence by [1, p.\ 330, Theorem 25.2, and p.\  383,
Theorem 29.4],
$$\cL\left(X^{(n)}_{k(1)},X^{(n)}_{k(2)},\dots,X^{(n)}_{k(p)}\right)
\Longrightarrow
\cL\left(X_{k(1)},X_{k(2)},\dots,X_{k(p)}\right)\,\,\hbox{as $n\to
\infty$} \eqno{(3.20)}$$ (weak convergence on $\R^p$). Recall the
elementary fact that for weak convergence on $\R^p$, the limiting
distribution is unique.

For a given $j\in\Z$ and a given $m\in\N$, one has that
$$\cL\left(X^{(n)}_{j+1},X^{(n)}_{j+2},\dots,X^{(n)}_{j+m}\right) =
\cL\left(X^{(n)}_1,X^{(n)}_2,\dots,X^{(n)}_m\right)$$
for each $n\ge 1$ by Lemma 3.6, and hence by (3.20) one has that
$$\cL\left(X_{j+1},X_{j+2},\dots,X_{j+m}\right) =
\cL\left(X_1,X_2,\dots,X_m\right).$$
Hence the sequence $X$ is
strictly stationary.

Since $\cL(X^{(n)}_0) = \lambda_{\unp3}$ for each $n\in\N$ by Lemma
3.6, it follows from (3.20) that $\cL(X_0) = \lambda_{\unp3}$. Thus
property (A) in Theorem 1.1 holds.

In the case $p=L-1$, the left side of (3.20) is
$\lambda^{[L-1]}_{\unp3}$ for each $n\in\N$ by Lemma 3.6 (see the
end of step 2.4(A)), and hence the right side of (3.20) is also
$\lambda^{[L-1]}_{\unp3}$. Thus property (B) (here
$(L-1)$--tuplewise independence) in Theorem 1.1 holds (see step
2.1).

Property (C) in Theorem 1.1 holds by a similar argument involving
Lemma 3.6 and eq.\ (3.20) with $p$ arbitrarily large.
\medskip

\noindent{\bf Step 3.11.} The final task is to prove property (D) in
Theorem 1.1. Recall the notations for partial sums in eq.\ (1.1).

By property (A) (proved above) in Theorem 1.1, and a routine
calculation, $EX_0 =0$, $EX^2_0 =1$, and $EX^4_0 =9/5$. Recall eq.\
(2.1) and the property of $(L-1)$--tuplewise independence (proved
above) for the sequence $X$. One has that
$$\eqalignno{
\forall\, n\ge 1,\quad &E\left(n^{-1/2}S(X,n)\right) =0
\quad\hbox{and} &(3.21) \cr \forall\, n\ge 1, \quad
&E\left(n^{-1/2}S(X,n)\right)^2 =1. &(3.22) \cr}$$ Further, by the
well known, elementary calculation in [1, p.\ 85, proof of Theorem
6.1] (which requires only 4--tuplewise independence), one has that
for each $n\in\N$,
$$\eqalign{
E(S(X,n))^4 &= n\cdot EX^4_0 +3n(n-1)(EX^2_0)^2 \cr
&= n\cdot (9/5) + (3n^2 -3n) \cdot 1 \le 3n^2. \cr}$$
Hence
$$\forall\, n\ge 1,\quad E\left( n^{-1/2}S(X,n)\right)^4 \le 3.
\eqno{(3.23)}$$
By (say) (3.22) and Chebyshev's inequality, the family of
distributions of the normalized partial sums $n^{-1/2}S(X,n)$,
$n\in\N$ is tight.

Now for the proof of property (D) in Theorem 1.1, suppose $Q$ is an
infinite subset of $\N$. Because of tightness, there exists an
infinite set $T\subset Q$ and a probability measure $\mu$ on $\R$
(both $T$ and $\mu$ henceforth fixed) such that
$$n^{-1/2}S(X,n) \Longrightarrow \mu \,\,\hbox{as $n\to
\infty$, $n\in T$.} \eqno{(3.24)}$$
To complete the proof of property (D) in Theorem 1.1, our task now is
to show that $\mu$ is neither degenerate nor normal.

Because of (3.23), (3.24), and [1, p.\ 338, the Corollary], one has
by (3.21) and (3.22) that
$$\int_{x\in\R} x\mu(dx) =0 \quad \hbox{and}\quad \int_{x\in\R}
x^2\mu(dx) =1. \eqno{(3.25)}$$ Hence the probability measure $\mu$
has positive variance and is therefore nondegenerate.

If $\mu$ were normal, then by (3.25) it would have to be the $N(0,1)$
distribution.  Because of (3.24) and [1, p.\ 334, Corollary 1, and
p.\ 338, Theorem 25.11], to show that $\mu$ fails to be normal, it
now suffices to show (see the sentence preceding Lemma 3.2) that
$$\forall\, n\ge 2L, \quad E\left(n^{-1/2}S(X,n)\right)^L \le EZ^L
-8^{-L/2} \cdot L! \cdot L^{-L}. \eqno{(3.26)}$$

Let $M\ge 2L$ be arbitrary but fixed. Let $m\in\N$ be the integer
such that
$$2L^m \le M <2L^{m+1}. \eqno{(3.27)}$$
By (3.27) and Lemma 3.7, for all $n\ge m$,
$$E\left[S(X^{(n)},M)\right]^L \le E\left[ M^{1/2}Z\right]^L -2^{-L}
\cdot L! \cdot L^{(m-1)L/2}. \eqno{(3.28)}$$

Now for each $n\in\N$, $[S(X^{(n)},M)]^L \le [3^{1/2}M]^L$ a.s.\
(since $|X^{(n)}_k| \le 3^{1/2}$ a.s.\ for each $k\in\Z$ by Lemma
3.6), and hence by (3.19) and dominated convergence,
$E[S(X^{(n)},M)]^L \to E[S(X,M)]^L$ as $n\to \infty$. Hence by
(3.28),
$$E[S(X,M)]^L \le E[M^{1/2}Z]^L - 2^{-L} \cdot L! \cdot
L^{(m-1)L/2}.$$ (That also follows simply from Fatou's Lemma.) Hence
by the second inequality in (3.27),
$$\eqalign{
E\left[M^{-1/2}S(X,M)\right]^L &\le EZ^L - 2^{-L} \cdot L! \cdot
L^{(m-1)L/2} \cdot M^{-L/2} \cr &\le EZ^L -2^{-L}\cdot L! \cdot
L^{(m-1)L/2} \cdot (2L^{m+1})^{-L/2} \cr &= EZ^L -2^{-3L/2}\cdot L!
\cdot L^{-L}. \cr}$$ Since $M\ge 2L$ was arbitrary, eq.\ (3.26)
holds. That completes the proof that $\mu$ is non-normal, and that
property (D) in Theorem 1.1 holds. The proof of Theorem 1.1 is
complete.
\bigskip
\bigskip
\noindent{\bf Acknowledgement.} The authors thank Jon Aaronson and
Benjamin Weiss for pointing out the reference [6], and Ron Peled for
pointing out the reference [8].

\vfill\eject
 \centerline{\bf References}
\medskip

\refs{[1]} P.\  Billingsley, {\it Probability and Measure}, third
ed., Wiley, New York, 1995.

\refs{[2]} R.C.\  Bradley, A stationary, pairwise independent,
absolutely regular sequence for which the central limit theorem
fails, {\it Probab.\  Th.\  Rel.\  Fields} 81 (1989), 1-10.

\refs{[3]} R.C.\ Bradley, On a stationary, triple-wise independent,
absolutely regular counterexample to the central limit theorem, {\it
Rocky Mountain J.\ Math.} 37 (2007), 25-44.

\refs{[4]} R.C.\ Bradley, A strictly stationary, ``causal,''
5--tuplewise independent counterexample to the central limit
theorem, (Revised version in preparation).

\refs{[5]} N.\ Etemadi, An elementary proof of the strong law of
large numbers, {\it Z.\ Wahrsch.\ verw.\ Gebiete} 55 (1981),
119-122.

\refs{[6]} L.\ Flaminio, Mixing $k$--fold independent processes of
zero entropy, {\it Proc. Amer.\ Math.\ Soc.} 118 (1993), 1263-1269.

\refs{[7]} S. Janson, Some pairwise independent sequences for which
the central limit theorem fails, {\it Stochastics} 23 (1988),
439-448.

\refs{[8]} F.S.\ MacWilliams and N.J.A.\ Sloane, {\it The Theory of
Error--Correcting Codes,} North Holland, Amsterdam, 1977.

\refs{[9]} A.R.\ Pruss, A bounded $N$-tuplewise independent and
identically distributed counterexample to the CLT, {\it Probab.\
Th.\ Rel.\ Fields} 111 (1998), 323-332.

\refs{[10]} A.R.\ Pruss, private communication to R.C.B., 2007.

\bye